\documentclass[12pt]{article}
\usepackage{a4,times,amsmath,amscd,amsfonts,amsthm,latexsym, amssymb, oldgerm}

\newtheorem{thm}{\qquad\sc Theorem} \newcommand{\bt}{\begin{thm}}\newcommand{\et}{\end{thm}}
\newtheorem{lemma}[thm]{\qquad \sc Lemma}\newcommand{\bl}{\begin{lemma}} \newcommand{\el}{\end{lemma}}
\newcommand{\pf}{{\qquad\it Proof:\quad}}
\newtheorem{prop}[thm]{\qquad\sc Proposition} \newcommand{\bp}{\begin{prop}}\newcommand{\ep}{\end{prop}}
\newtheorem{cor}[thm]{\qquad \sc Corollary} \newcommand{\bc}{\begin{cor}}\newcommand{\ec}{\end{cor}}
\newtheorem{rem}[thm]{\qquad\it Remark} \newcommand{\br}{\begin{rem}}\newcommand{\er}{\end{rem}}
\newtheorem{con}[thm]{\qquad\sc Conjecture} \newcommand{\bcon}{\begin{con}}\newcommand{\econ}{\end{con}}

\begin{document}

\begin{center}
{\bf SECOND ORDER AVERAGE ESTIMATES ON LOCAL DATA OF CUSP FORMS}
\end{center}

\begin{center}
{\large {\sc Farrell Brumley}}\footnote{The author was supported by a postdoctoral fellowship from the CNRS.}
\end{center}

\bigskip

\begin{abstract}
We specify sufficient conditions for the square modulus of the local parameters of a family of ${\rm GL}_n$ cusp forms to be bounded on average.  These conditions are global in nature and are satisfied for $n\leq 4$.  As an application, we show that Rankin-Selberg $L$-functions on ${\rm GL}_{n_1}\times {\rm GL}_{n_2}$, for $n_i\leq 4$, satisfy the standard convexity bound.
\end{abstract}

\begin{center}
\section{\rm Introduction}
\end{center}

\qquad Let $F$ be a number field and $\mathbb{A}$ its ring of adeles.  Let $\pi=\otimes_v\pi_v$ be a cuspidal automorphic representation of ${\rm GL}_n(\mathbb{A})$.  At each finite place $v=\textfrak{p}$ of $F$ there is associated with $\pi_\textfrak{p}$ a semisimple conjugacy class $A_\pi (\textfrak{p})$ in ${\rm GL}_n(\mathbb{C})$, the matrix of local (Langlands) parameters $A_\pi (\textfrak{p})={\rm diag} (\alpha_\pi(\textfrak{p},1),\ldots ,\alpha_\pi (\textfrak{p},n))$.  The Ramanujan conjecture states that $|\alpha_\pi (\textfrak{p},i)|\leq 1$ for all $1\leq i\leq n$, with strict equality when $\pi_\textfrak{p}$ is unramified.

\qquad One may use the information supplied by the Ramanujan Conjecture to derive important analytic results for $L$-functions.  In doing so one trades strong pointwise information for results that often have more to do with the average behavior of the local parameters, with $\alpha_\pi(\textfrak{p},i)$ ranging over primes $\textfrak{p}$ and possibly over $\pi$ is some family.  One such consequence of the Ramanujan Conjecture is that any (appropriately normalized) $L$-function associated to $\pi$ satisfies the optimal estimate $O_\epsilon (C(\pi)^\epsilon)$ on ${\rm Re}(s)>1$.  The quantity $C(\pi)$ is the analytic conductor of $\pi$ (see Section \ref{global} for the definition).  When this property holds, we say that the $L$-function satisfies the standard convexity bound.

\qquad One technique used to demonstrate optimal bounds on sums of positive coefficients was introduced by Iwaniec [6] for cusp forms $\pi$ on ${\rm GL}_2$.  Iwaniec uses a linearization process to show that if the coefficients $\lambda (\textfrak{n},\pi)$ of $L(s,\pi)$ don't begin to show $O(1)$ behavior by the time ${\rm N}\textfrak{n}$ is of size $O_\epsilon (C(\pi)^\epsilon)$ then this late excess will so propogate through the remaining coefficients via their multiplicative relations as to contradict the polynomial control granted by the Rankin-Selberg theory.  Molteni [16], working out the difficult combinatorics involved in implementing Iwaniec's idea in full generality, was then able to show that for $\pi$ any cusp form on ${\rm GL}_n$ the principal $L$-function $L(s,\pi)$ satisfies the standard convexity bound.

\qquad To apply the same reasoning to the Rankin-Selberg $L$-function $L(s,\pi\times\pi)$ requires a more delicate analysis.  It may come as a surprise to some that despite the recent breakthroughs in certain cases of {\it sub}convexity, it is still not known in complete generality and under no assumptions that $L(s,\pi\times\pi)$ satisfies the standard convexity bound.  Molteni [16] went some way toward this goal by showing that for $\pi$ any cusp form on ${\rm GL}_n$ as long as 

\begin{equation}\label{molt}
|\alpha_\pi (\textfrak{p},i)|\ll {\rm N}\textfrak{p}^{1/4} 
\end{equation}

for all but finitely many primes $\textfrak{p}$ and $1\leq i\leq n$ then $L(s,\pi\times\pi)$ satisfies the standard convexity bound (see his Hypothesis (R$'$)).  At present, however, bounds of this quality are known only for cusp forms on ${\rm GL}_2(\mathbb{A})$ where we have $|\alpha_\pi (\textfrak{p},i)|\ll {\rm N}\textfrak{p}^{1/9}$ [11].  

\qquad In this paper we remove hypothesis (\ref{molt}) in certain cases, proving that $L(s,\pi_1\times\pi_2)$ satisfies the standard convexity bound for pairs $(\pi_1,\pi_2)$ on ${\rm GL}_{n_1}\times {\rm GL}_{n_2}$ for $n_i\leq 4$.  This improvement upon the range given by Molteni's work is due to a greater emphasis on global infomation and benefits from some recent advances in functoriality.  Throughout the paper we take pains to describe what happens on higher rank in an effort to compare the strengths of our method with those of other approaches.

\bigskip

\subsection{\rm Main Theorem}

\qquad We consider a Dirichlet series which acts as a majorizer of $L(s,\pi\times\pi)$.  For $\pi$ any cuspidal representation of ${\rm GL}_n(\mathbb{A})$, define

\begin{equation*}
L(s,\pi,|{\rm max} |^2):=\sum_\textfrak{n}\lambda (\textfrak{n},\pi,|{\rm max}|^2){\rm N}\textfrak{n}^{-s}:=\prod_\textfrak{p}\sum_{r\geq 0}\underset{i}{\rm max}|\alpha_\pi (\textfrak{p},i)|^{2r}\ {\rm N}\textfrak{p}^{-rs}.
\end{equation*}

We shall specify sufficient conditions under which this $L(s,\pi,|{\rm max}|^2)$ is $O_\epsilon (C(\pi)^\epsilon)$ on ${\rm Re}(s)>1$.  We call this estimate the convexity bound {\it at $s=1$}, detailing the specific point in this case since $L(s,\pi,|{\rm max}|^2)$, lacking a functional equation, does not allow for an interpolation to points to the left of 1.

\qquad This function $L(s,\pi,|{\rm max}|^2)$ has the advantage over $L(s,\pi\times\tilde \pi)$ of being completely multiplicative in its coefficients.  As we shall see in Proposition \ref{barcelona}, Dirichlet series whose coefficients are positive and completely multiplicative can be subjected to Iwaniec's boot-strapping method with no additional assumption on the size of their coefficients.

\qquad The disadvantage of working with $L(s,\pi,|{\rm max}|^2)$ is that it is not an $L$-function coming from an automorphic form, making its analytic properties hard to unearth.  To remedy this problem, we majorize $\lambda (\textfrak{p},\pi,|{\rm max}|^2)$ for $\textfrak{p}$ at which $\pi_\textfrak{p}$ is unramified by a sum of the absolute values of certain more naturally arising coefficients (see Proposition \ref{point}).  In doing so, we make use of the fact that for unramified $\textfrak{p}$ the unitarity of $\pi_\textfrak{p}$ restricts the number of roots that can possibly violate the Ramanujan conjecture.  The matrix of Satake parameters $A_\pi (\textfrak{p})$ in this case is forced to lie in the same semi-simple conjugacy class as $\overline{A_\pi (\textfrak{p})^{-1}}$, meaning that only $\lfloor n/2\rfloor$ of the $n$ roots can have size greater than 1.  The function $\lfloor\ \cdot\ \rfloor$ is the ``floor'' function outputting the largest integer less than or equal to the imput value.

\qquad The following is our main theorem.  For the definition of a {\it strong} isobaric lift consult Section \ref{unit}.

\bt\label{olivia} Let $\pi$ a cuspidal representation of ${\rm GL}_n(\mathbb{A})$.  For an integer $j\geq 2$ denote by $\land^j$ the exterior $j$-power representation of ${\rm GL}_n(\mathbb{C})$.  Assume that for every $2\leq j\leq\lfloor n/2\rfloor$ there exists a strong $\land^j$-isobaric lift.  Then $L(s,\pi,|{\rm max}|^2)$ satifies the convexity bound at $s=1$.\et

\subsection{\rm Applications}

\qquad When $n=2$ or $3$, the conditions of Theorem \ref{olivia} are empty, so the conclusion automatically holds.  When $n=4$ or $5$ the sole condition in is that there exists a strong isobaric $\land^2$ lift.  For $n= 4$, this condition was proven by Kim [9, Proposition 5.3.1], and so we can state unconditionally the following practical result.  For cusp forms on ${\rm GL}_3$ and ${\rm GL}_4$, this result is new.

\bc\label{nice} Let $\pi_i$ be cuspidal representations of ${\rm GL}_{n_i}(\mathbb{A})$ where $n_i\leq 4$ for $i=1, 2$.  Then $L(s,\pi_1\times\pi_2)$, as well as $L(s,\pi_i,\land^2)$ and $L(s,\pi_i,{\rm sym}^2)$ for $i=1, 2$, satisfy the standard convexity bound.\ec

\qquad The proof of Corollary \ref{nice}, essentially an application of the Cauchy-Schwartz inequality, is provided in sub-section \ref{cor-proof}.

\qquad We now given several examples where Theorem \ref{olivia} can be used to replace the Ramanujan Conjecture or the hypothetical bounds (\ref{molt}).  The first is a large sieve inequality for long sums of Fourier coefficients of cusp forms on ${\rm GL}_n$.  For cusp forms $\pi$ on ${\rm GL}_2/\mathbb{Q}$ and their images ${\rm sym}^2\pi$ under the Gelbart-Jacquet lift this large sieve inequality is a theorem, by Duke and Kowalski [2] in the level aspect when $\pi$ is holomorphic, by Luo [13] in the eigenvalue aspect when $\pi$ is a Maass form.  

\qquad For a parameter $Q\geq 1$ let $S_n(\leq Q)$ be the set of all cusp forms on ${\rm GL}_n/\mathbb{Q}$ with analytic conductor bounded by $Q$.  Under a remaining assumption giving polynomial growth on $S_n(\leq Q)$, the results of [2] and [13] can be extended to $n\leq 4$ using Theorem \ref{olivia}.

\bc\label{large sieve} Let $n=3$ or $4$.  Assume that there exists a number $d>0$ such that $|S_n(\leq Q)| = O(Q^d)$.  Let $\alpha=1-(n^2+1)^{-1}$ and for $\pi_1, \pi_2\in S_n(\leq Q)$ define $B=B(n)>0$ to be the exponent appearing in the convexity bound of $L(s,\pi_1\times\pi_2)$ at $s=\alpha$, so that $L(\alpha,\pi_1\times\pi_2)\ll_\epsilon Q^{B+\epsilon}$.  Then for any $\epsilon > 0$ the inequality

\begin{equation*}
\sum_{\pi\in S(\leq Q)}\bigg|\sum_{n\leq N} a_n \lambda (n,\pi)\bigg|^2\ll_\epsilon (NQ)^\epsilon (N+Q^{B+d}N^\alpha)\sum_{n\leq N} |a_n|^2 
\end{equation*}

holds for all complex numbers $(a_n)_{1\leq n\leq N}$.\ec

\qquad The proof of Corollary \ref{large sieve} is by a well-known duality argument.  We shall only sketch the details which pertain to the role of Theorem \ref{olivia}.  The two terms on the right-hand side of the large sieve inequality come from a majorization of an integral involving the Rankin-Selberg $L$-function at the point $s=1$ and along the line $s=\alpha+\epsilon$ for small $\epsilon>0$.  At $s=1$ one uses Theorem \ref{olivia} to show that the residue of $L(s,\pi\times\tilde\pi)$ is $\ll_\epsilon Q^\epsilon$.  On the line $s=\alpha+\epsilon$ Theorem \ref{olivia} is used to establish the convergence and negligibility of the correction factor that relates the bilinear Rankin-Selberg $L$-function to the true convolution.  This correction factor, labeled $H(s,\pi_1\times\pi_2)$ in [2], is the product $\prod_p H_p(s,\pi_1\times\pi_2)$ of polynomials in $p^{-s}$ whose coefficients are symmetric polynomials in the local roots of $\pi_1$ and $\pi_2$ and whose linear term is zero.  The product converges on ${\rm Re}(s)>\alpha$ and satisfies $H(s,\pi_1\times\pi_2)\ll_\epsilon Q^\epsilon$ in this region by the Luo-Rudnick-Sarnak bounds (see display (\ref{LRS bounds})) and Theorem \ref{olivia}.

\qquad As a final application we state a zero-density theorem for principal $L$-functions of cusp forms on ${\rm GL}_n$ for $n=3$ or $4$.  Kowalski and Michel [10] have proven a quite powerful zero-density statement near the line 1 for $L$-functions of general ${\rm GL}_n$ cusp forms over $\mathbb{Q}$.  They assume the bounds (\ref{molt}) on the local parameters, in this case to prove a mean-value theorem with pseudo-characters {\it \`a la Selberg}.  As in Corollary \ref{large sieve} we may remove this assumption for $n=3$ and $4$ using our Theorem \ref{olivia}.  See the introduction to [10] for applications of this type of zero-density statement to moments of $L(1,\pi)$ for certain families of $\pi$.

\qquad Let $M(\alpha, T) = \{ z\in \mathbb{C} |\ {\rm  Re}(z)\geq\alpha\ \text{and}\ | {\rm Im}(z)|\leq T\}$ for $\alpha\in\mathbb{R}$, $T\geq 0$. For any cuspidal automorphic representation $\pi$ of ${\rm GL}_n/\mathbb{Q}$, we let 

\begin{equation*}
N(\pi; \alpha, T) = |\{\rho\in M(\alpha, T) | L(\pi, \rho) = 0\}| 
\end{equation*}

(zeros counted with multiplicity).

\bc Let $n=3$ or $4$.  Let $S_n(\leq Q)$, $Q\geq 1$, be as above.  Assume that there exists a number $d>0$ such that $|S_n(\leq Q)|=O(Q^d)$.  Let $\alpha>1-(n^2+1)^{-1}$ and $T\geq 2$.  Then there exist constants $c, B > 0$, depending only on $n$ and $d$, such that

\begin{equation*}
\sum_{\pi\in S_n(\leq Q)} N(\pi;\alpha , T)\ll T^B Q^{c(1-\alpha)}
\end{equation*}

for all $Q\geq 1$.  The implied constant depends only on the choice of $c$.\ec

\subsection{\rm Strength of method}

\qquad We have tried in this paper to give the reader an idea of the strengths of our method relative to other approaches.  It is for this reason that, despite the extremity of its hypotheses, Theorem \ref{olivia} was stated for general $n$.

\qquad One weakness of the method we outline is that there is much information loss in passing from the conclusion of Theorem \ref{olivia} to Corollary \ref{nice}.  This lost information is hard to quantify, and it is not at all clear that additional applications could be gleaned from the stronger result.  To see the information loss, imagine trying to reverse the logic to deduce the convexity bound at $s=1$ of $L(s,\pi,|{\rm max}|^2)$ from that of $L(s,\pi,{\rm sym}^2)$ and $L(s,\pi,\land^2)$ alone.  Note that convexity for these latter two implies the same for their product $L(s,\pi\times\pi)$, and indeed for any Rankin-Selberg pair $L(s,\pi_1\times\pi_2)$ by the Cauchy-Schwartz inequality.  But as $n$ gets large ($n$ greater that 5 will work), our Proposition \ref{point} shows that many more representations are needed to control the modulus-squared of the roots.   

\qquad Theorem \ref{olivia} thus seems best suited for $n\leq 5$ where the degree of the representations whose automorphy we assume is no larger than those to whose $L$-functions we apply the result.  But even for $n\leq 5$, where the standard and the exterior square representation suffice to control the square-modulus of the roots, there is information loss simply by the reduction to a completely multiplicative Dirichlet series.  When applying Theorem \ref{olivia} to $L(s,\pi\times\tilde\pi)$ for example, we are using the quantity $(\lambda(\textfrak{p},\pi\times\tilde\pi)+|\lambda (\textfrak{p},\pi,\land^2)|)^r$ to control $\lambda (\textfrak{p}^r,\pi\times\tilde\pi)$ for every $r\geq 0$.  When $r=1$ the presence of the exterior square is clearly unnecessary.  By treating all coefficients with essentially the same majorization, we neutralize the otherwise helpful effect of interior cancellation among the roots that might lead to some coefficients being small or zero.

\qquad Let us say more about other methods for proving the standard convexity bound for Rankin-Selberg $L$-functions.  The most direct way to force the convexity bound for $L(s,\pi\times\pi)$ is by assuming the existence of both an exterior and symmetric square lift.  For it is clear from the identity $L(s,\pi\times\pi)=L(s,\pi,{\rm sym}^2)L(s,\pi,\land^2)$ that the convexity bound for $L(s,\pi\times\pi)$ follows from that of both $L(s,\pi,{\rm sym}^2)$ and $L(s,\pi,\land^2)$; by the results in [16] again, this would follow from the (isobaric) automorphy of both ${\rm sym}^2\pi$ and $\land^2\pi$.  Hence $L(s,\pi\times\pi)$, where $\pi$ is a cusp form on ${\rm GL}_2(\mathbb{A})$, satisfies the standard convexity bound for yet another reason: the Gelbart-Jacquet lift [3].

\qquad In Section \ref{unit} we are able to shed some light on the relation between the hypothetical bounds (\ref{molt}) on the local roots and the direct assumption of functoriality of both $\land^2$ and ${\rm sym}^2$.  In Corollary \ref{over with} we show that for $n\leq 5$ the assumption of both functorial lifts is stronger than the condition (\ref{molt}).  For $n>5$ no such implication can be made by our method.  In fact, working locally one unramified prime at a time, and using only unitarity as input, we show that many more functorial lifts are needed to break the $1/4$ exponent that Molteni requires.  Of course we lose lots of global information in setting up this implication, but it is interesting nonetheless to consider whether, for higher rank general linear groups, the existence of both functorial lifts, already such an extreme hypothesis, might actually be weaker than the bounds in (\ref{molt}).

\bigskip

\qquad {\it Acknowledgements.}  This paper represents a chapter in my doctoral dissertation.  I would like to thank my thesis advisor Peter Sarnak for suggesting this problem to me and Ramin Takloo-Bighash, Philippe Michel, and Akshay Venkatesh for their encouragement.

\begin{center}
\section{\rm Consequences of unitarity}\label{unit}
\end{center}

\qquad In the following proposition, we have chosen the exterior power lifts for simplicity.  The proof uses only properties on the size of the eigenvalues $\alpha_i$.  Alternative sets of representations of ${\rm GL}_n(\mathbb{C})$ may be chosen, though one would then have to take into consideration the arguments of the $\alpha_i$.

\bigskip

\bp\label{point} Let $n\geq 1$ and $m=\lfloor n/2\rfloor$.  There exists a constant $c_n>0$ depending only on $n$ such that for any matrix $A={\rm diag}(\alpha_1,\ldots ,\alpha_n)\in {\rm GL}_n(\mathbb{C})$ with $A^{-1}$ and $\overline{A}$ lying in the same semi-simple conjugacy class

\begin{equation}\label{prop display}
\underset{i}{\rm max}|\alpha_i|^2 \leq c_n\left( 1+\sum_{j=2}^m|{\rm Trace}(\land^j A)|^{2/j}\right).
\end{equation}

\ep

\pf The assumption on $A$ means that there is some permutation $\sigma$ of the indices such that $\alpha_i\overline{\alpha_{\sigma(i)}}=1$ for all $i$.  The elements may be ordered by their size, say

\begin{equation}\label{ordering}
|\alpha_1|\geq \cdots \geq |\alpha_m|\geq 1\geq \cdots \geq |\alpha_n|.
\end{equation}

We note that

\begin{equation}\label{trace}
{\rm Trace}(\land^j A)=\underset{1\leq i_1<\cdots <i_j\leq n}{\sum}\alpha_{i_1}\cdots\alpha_{i_j}.
\end{equation}

\qquad For the moment let $R_1,\ldots ,R_{m+1}$ be any array of positive real numbers satisfying $R_1=1$ and $0< R_i< 1$ for $2\leq i\leq m+1$.  It is clear that either

\bigskip

\qquad (i) there exists some $j\in\{ 1,\ldots m\}$ such that $|\alpha_i|\geq R_i|\alpha_1|$ for all $1\leq i\leq j$ and $|\alpha_{j+1}|\leq R_{j+1}|\alpha_1|$; or else

\bigskip

\qquad (ii) $|\alpha_{m+1}|\geq R_{m+1}|\alpha_1|$. 

\bigskip

In case (ii), we have $|\alpha_1|^2\leq R_{m+1}^{-2}|\alpha_m|^2\leq R_{m+1}^{-2}$ by (\ref{ordering}).

\qquad Now let $j$ be as in case (i).  The leading term in (\ref{trace}) is $\alpha_1\cdots\alpha_j$ which has size

\begin{equation*}
|\alpha_1\cdots\alpha_j|\geq \left(\prod_{i=1}^jR_i\right)|\alpha_1|^j.
\end{equation*}

From (\ref{ordering}) all other terms are bounded in absolute value by

\begin{equation*}
|\alpha_1\cdots\alpha_{j-1}\alpha_{j+1}|\leq R_{j+1}|\alpha_1|^j.
\end{equation*}

Thus we have

\begin{equation}\label{error}
{\rm Trace}(\land^j A)=\alpha_1\cdots\alpha_j+O((r_j-1)R_{j+1}|\alpha_1|^j),
\end{equation}

where the implied constant is bounded by 1 and $r_j$ be the number of terms present in ${\rm Trace}(\land^j A)$, that is $r_j=\# \{ (i_1,\cdots ,i_j)\ |\ 1\leq i_1<\cdots <i_j\leq n\}$.  The numbers $R_i$ should now be chosen to make the main term in (\ref{error}) dominate the error term.  For $1\leq j\leq m+1$, set

\begin{equation*}
R_j=\left( \prod_{i=1}^{j-1} R_i\right) r_{j-1}^{-1}.
\end{equation*}

Here we have put $r_0=1$.  Thus $R_1=1$, $R_2=r_1^{-1}$, $R_3=r_1^{-1}r_2^{-1}$, and so on.   By (\ref{error}) this implies

\begin{equation*}
|\alpha_1|^j\leq \left(\prod_{i=1}^j R_i\right)^{-1}\left(1-\frac{r_j-1}{r_j}\right)^{-1}|{\rm Trace}(\land^j A)|= R_{j+1}^{-1} |{\rm Trace}(\land^j A)|.
\end{equation*}

\qquad Since $R_1>R_2>\ldots >R_{m+1}$, we may encompass all cases by taking $c_n=R_{m+1}^{-2}$.  This completes the proof.\qed

\bigskip

\qquad Let $\pi_v$ be any irreducible admissible representation of ${\rm GL}_n (F_v)$.  Let $\phi_v: W'_{F_v}\to {\rm GL}_n(\mathbb{C})$ the parametrization of $\pi_v$ given by the local Langlands correspondence ([4], [5], [12]), where $W'_{F_v}$ is the Weil-Deligne group.  The composition $\phi_v\cdot\rho$ is then the parametrization of an irreducible admissible representation $\rho (\pi_v)$ of ${\rm GL}_N(F_v)$.  We may form the tensor product $\rho(\pi)=\otimes_v\rho(\pi_v)$ over all places $v$ of $F$.  The result is an irreducible admissible representation of ${\rm GL}_N(\mathbb{A})$. An automorphic representation $\Pi$ of ${\rm GL}_N(\mathbb{A})$ is called {\it isobaric} if $\Pi={\rm Ind}\ \sigma_1\otimes\cdots\otimes\sigma_k$, for cuspidal representations $\sigma_i$ of ${\rm GL}_{n_i}(\mathbb{A})$, where $n_1+\cdots +n_k=N$.  Langlands functoriality predicts that $\rho(\pi)$ is automorphic.

\qquad We shall call an automorphic (respectively, isobaric) representation $\Pi^\rho=\otimes_v\Pi^\rho_v$ of ${\rm GL}_N(\mathbb{A})$ a {\it weak $\rho$-automorphic (respectively, -isobaric) lift} of a cuspidal representation $\pi=\otimes_v\pi_v$ of ${\rm GL}_n(\mathbb{A})$ if there exists a finite set $S_\pi$ of places, including the finite places $v$ at which $\pi_v$ is ramified, such that $\Pi^\rho_v\simeq\rho (\pi_v)$ for all $v\notin S_\pi$.  The lift is said to be {\it strong} if $S_\pi$ can be taken independent of $\pi$.

\qquad Using our Proposition \ref{point} and additional assumptions on the existence of certain weak isobaric lifts, the Luo-Rudnick-Sarnak [14] bounds 

\begin{equation}\label{LRS bounds}
|\alpha_\pi (\textfrak{p},i)|\leq {\rm N}\textfrak{p}^{1/2-(n^2+1)^{-1}},
\end{equation}

valid for all primes $\textfrak{p}$ and $1\leq i\leq n$, can be dramatically improved.

\bc\label{over with} Let $\pi$ be a cuspidal automorphic representation of ${\rm GL}_n(\mathbb{A})$.  Put $m=\lfloor n/2\rfloor$.  Assume that there exists a weak ${\rm sym}^2$-isobaric lift and $\land^j$-isobaric lifts of $\pi$ for all $2\leq j\leq m$.  There exists a $\delta_n>0$ such that

\begin{equation*}
{\rm N}\textfrak{p}^{-1/4+\delta_n}\ll |\alpha_\pi (\textfrak{p},i)|\ll {\rm N}\textfrak{p}^{1/4-\delta_n}
\end{equation*}

for $1\leq i\leq n$ and almost all primes $\textfrak{p}$.
\ec

\pf Using the trivial inequality $|{\rm Trace}(A)|^2\leq |{\rm Trace}({\rm sym}^2 A)|+|{\rm Trace}(\land^2 A)|$, Proposition \ref{point} gives

\begin{equation}\label{stop}
\underset{i}{\rm max}|\alpha_\pi (\textfrak{p},i)|^2\ll |{\rm Trace}({\rm sym}^2 A_\pi(\textfrak{p}))|+\sum_{j=2}^m|{\rm Trace}(\land^j A_\pi (\textfrak{p}))|
\end{equation}

for all $\textfrak{p}$ such that $\pi_\textfrak{p}$ is unramified.  Let $\rho: {\rm GL}_n(\mathbb{C})\to {\rm GL}_d(\mathbb{C})$ be a polynomial representation and $\Pi=\Pi(\rho)$ a weak $\rho$-isobaric lift of $\pi$.  Let $\Pi={\rm Ind}\ \sigma_1\otimes\cdots\otimes \sigma_k$ with $\sigma_i$ a cusp form on ${\rm GL}_{d_i}(\mathbb{A})$ and $d_1+\cdots +d_k=d$.  We have

\begin{equation*}
|{\rm Trace}\rho (A_\pi (\textfrak{p}))|=\bigg|\sum_{i=1}^k {\rm Trace}( A_{\sigma_i}(\textfrak{p}))\bigg|\leq \sum_{i=1}^k d_i {\rm N}\textfrak{p}^{1/2-(d_i^2+1)^{-1}}\leq d{\rm N}\textfrak{p}^{1/2-(d^2+1)^{-1}},
\end{equation*}

for all $\textfrak{p}\notin S_\pi$ by the Luo-Rudnick-Sarnak bounds (\ref{LRS bounds}).  With $d={\rm max}\{ {\rm deg}({\rm sym}^2), {\rm deg}(\land^m)\}$ we apply this upper bound to each summand on the right hand side of (\ref{stop}) to get

\begin{equation*}
\underset{i}{\rm max}|\alpha_\pi (\textfrak{p},i)|^2\ll {\rm N}\textfrak{p}^{1/2-(d^2+1)^{-1}}.
\end{equation*}

Corollary \ref{over with} then follows from the unitarity of $\pi_\textfrak{p}$.  \qed

\begin{center}
\section{\rm Global estimates}\label{global}
\end{center}

\qquad We construct a Dirichlet series which will be the focus of our attention for the rest of this paper.  Let $\pi$ be a cuspidal representation of ${\rm GL}_n(\mathbb{A})$.  Define

\begin{equation*}
L(s,\pi,|{\rm max} |^2):=\sum_\textfrak{n}\lambda (\textfrak{n},\pi,|{\rm max}|^2){\rm N}\textfrak{n}^{-s}:=\prod_\textfrak{p}\sum_{r\geq 0}\underset{i}{\rm max}|\alpha_\pi (\textfrak{p},i)|^{2r}\ {\rm N}\textfrak{p}^{-rs}.
\end{equation*}

Let $T$ be a finite set of primes such that $\textfrak{p}\notin T$ implies $\pi_\textfrak{p}$ is unramified.  Denote by $\textfrak{t}$ the (square-free) ideal which is the product of all primes in $T$.  Write

\begin{equation}\label{unramified L}
L_T(s,\pi,|{\rm max}|^2)=\prod_{\textfrak{p}\nmid\textfrak{t}}\sum_{r\geq 0}\underset{i}{\rm max}|\alpha_\pi(\textfrak{p},i)|^2{\rm N}\textfrak{p}^{-rs}=\sum_{(\textfrak{n},\textfrak{t})=1}\lambda(\textfrak{n},\pi,|{\rm max}|^2){\rm N}\textfrak{n}^{-s},
\end{equation}

and

\begin{equation}\label{ramified L}
L^T(s,\pi,|{\rm max}|^2)=\prod_{\textfrak{p}|\textfrak{t}}\sum_{r\geq 0}\underset{i}{\rm max}|\alpha_\pi(\textfrak{p},i)|^2{\rm N}\textfrak{p}^{-rs}.
\end{equation}

The following proposition is a consequence of Proposition \ref{point} and the Luo-Rudnick-Sarnak bounds (\ref{LRS bounds}).  The full strength of the bounds (\ref{LRS bounds}) is actually not used until the calculations involving ramified primes in Proposition \ref{catherine}.

\bigskip

\bp\label{harken back} Let $\pi$ be a cuspidal representation of ${\rm GL}_n(\mathbb{A})$.   Then

\begin{equation*}
L_T(s,\pi,|{\rm max}|^2)\ll \sum_{(\textfrak{n},\textfrak{t})=1}^{\qquad\flat}{\rm N}\textfrak{n}^{-\sigma}\sum_{(\textfrak{n},\textfrak{t})=1}^{\qquad\flat} \lambda(\textfrak{n},\pi\times\tilde\pi) {\rm N}\textfrak{n}^{-\sigma} \prod_{j=2}^m\sum_{(\textfrak{n},\textfrak{t})=1}^{\qquad\flat}|\lambda (\textfrak{n},\pi,\land^j)|{\rm N}\textfrak{n}^{-\sigma}
\end{equation*}

uniformly on ${\rm Re}(s)=\sigma\geq\sigma_0>1$.  The $\flat$ sign in the above sums indicates a restriction to square-free integral ideals.
\ep

\pf From the $\textfrak{p}$-th factor of $L_T(s,\pi,|{\rm max}|^2)$ we may extract a linear term to obtain

\begin{equation*}
\sum_{r\geq 0}\lambda(\textfrak{p}^r,\pi,|{\rm max}|^2){\rm N}\textfrak{p}^{-rs}=(1+\lambda(\textfrak{p},\pi,|{\rm max}|^2){\rm N}\textfrak{p}^{-s})\sum_{r\geq 0}\lambda (\textfrak{p}^{2r},\pi,|{\rm max}|^2){\rm N}\textfrak{p}^{-2rs}.
\end{equation*}

The bounds $\lambda (\textfrak{p},\pi,|{\rm max}|^2)\leq {\rm N}\textfrak{p}$, guaranteed to hold by (\ref{LRS bounds}), are strong enough to show convergence of the above geometric series to the right of 1.  (We shall need the full strength of the bounds (\ref{LRS bounds}) to treat $L^T(s,\pi,|{\rm max}|^2)$ in Proposition \ref{catherine}.)  Thus

\begin{equation*}
L_T(s,\pi,|{\rm max} |^2)\ll \prod_{\textfrak{p}\nmid\textfrak{t}}(1+\lambda(\textfrak{p},\pi,|{\rm max}|^2){\rm N}\textfrak{p}^{-\sigma})
\end{equation*}

uniformly on ${\rm Re}(s)=\sigma\geq\sigma_0>1$.  By Proposition \ref{point},

\begin{equation*}
\lambda(\textfrak{p},\pi,|{\rm max}|^2)\ll 1+\lambda(\textfrak{p},\pi\times\tilde\pi)+\sum_{j=2}^m|\lambda (\textfrak{p}, \pi,\land^j)| \qquad\text{for}\qquad \textfrak{p}\nmid \textfrak{t}.
\end{equation*}

Applying this majorization gives

\begin{align*}
&1+\lambda(\textfrak{p},\pi,|{\rm max}|^2){\rm N}\textfrak{p}^{-\sigma} \ll 1+ {\rm N}\textfrak{p}^{-\sigma}+\lambda(\textfrak{p},\pi\times\tilde\pi){\rm N}\textfrak{p}^{-\sigma}+\sum_{j=2}^m|\lambda (\textfrak{p},\pi,\land^j)|{\rm N}\textfrak{p}^{-\sigma}\\
& \qquad\qquad\qquad \leq (1+{\rm N}\textfrak{p}^{-\sigma})(1+ \lambda(\textfrak{p},\pi\times\tilde\pi){\rm N}\textfrak{p}^{-\sigma})\prod_{2\leq j\leq m}(1+ |\lambda (\textfrak{p},\pi,\land^j)|{\rm N}\textfrak{p}^{-\sigma}).
\end{align*}

Taking the product over all $\textfrak{p}\nmid\textfrak{t}$ we obtain the proposition.\qed 

\bc\label{inutile} Let $\pi$ be a cuspidal representation of ${\rm GL}_n(\mathbb{A})$.  If $L(s,\pi,\land^2)$ converges absolutely to the right of 1 for all $2\leq j\leq \lfloor n/2\rfloor$, then $L(s,\pi,|{\rm max}|^2)$ converges to the right of 1. \ec

\qquad When $n=2$ or $3$, the conditions of Corollary \ref{inutile} are empty, so the conclusion automatically holds.  When $n=4$ or $5$ the sole condition is that $L(s,\pi,\land^2)$ be absolutely convergent to the right of 1.  For $n= 4$, this property is proven by Kim in [9, Proposition 6.2].  Thus $L(s,\pi,|{\rm max}|^2)$ converges to the right of 1 for any $\pi$ on ${\rm GL}_n$ for $n\leq 4$.

\qquad When applied to Dirichlet series which arise naturally in the theory of automorphic forms, Corollary \ref{inutile} gives no new information.  For to deduce the absolute convergence to the right of 1 of $L(s,\pi,\land^2)$ from that of $L(s,\pi,|{\rm max}|^2)$ is just to repeat one of the hypotheses from which we derived the latter fact.  The same can be said for $L(s,\pi,{\rm sym}^2)$.  At this point, the loss of information in passing from $L(s,\pi,|{\rm max}|^2)$ to either $L(s,\pi,\land^2)$ or $L(s,\pi,{\rm sym}^2)$ is just too great.  The true strength of Proposition \ref{harken back} will presently be seen to lie in questions regarding uniformity in the analytic conductor of $\pi$.  

\subsection{\rm Gaining uniformity}

\qquad Denote the local parameters of $\pi$ at the infinite place $v$ by $\mu_\pi(v,i)$, $1\leq i\leq n$.  Let $q(\pi)$ be the conductor of $\pi$ and define the analytic conductor to be $C(\pi)=q(\pi)\lambda_\infty(\pi)$ where

\begin{equation*}
\lambda_\infty(\pi)=\prod_{v=\infty}\prod_{i=1}^n (1+|\mu_\pi (v,i)|).
\end{equation*}

For a pair of cusp forms $\pi_1, \pi_2$ on ${\rm GL}_{n_1}(\mathbb{A})$ and ${\rm GL}_{n_2}(\mathbb{A})$ we define the analytic conductor using the parameters at infinity present in the gamma factors of the completed $L$-function.  That is, for an infiinite place $v$,

\begin{equation*}
L_v(s,\pi_{1,v}\times\pi_{2,v})=\prod_{i=1}^{n_1}\prod_{j=1}^{n_2}\Gamma_{F_v}(s+\mu_{\pi_1\times\pi_2}(v,i,j))
\end{equation*}

for complex numbers $\mu_{\pi_1\times\pi_2}(v,i,j)$.  Above we have used the standard notation $\Gamma_{\mathbb{R}}(s)=\pi^{-s/2}\Gamma (s/2)$ and $\Gamma_{\mathbb C} (s)=2(2\pi)^{-s}\Gamma (s)$.  When the infinite place $v$ is unramified for either $\pi$ and $\pi'$ we have $\{\mu_{\pi\times\pi'}(v,i,j)\}=\{\mu_\pi(v,i)+\mu_{\pi'}(v,j)\}$.  We define the analytic conductor of $L(s,\pi_1\times\pi_2)$ to be $C(\pi_1\times\pi_2)=q(\pi_1\times\pi_2)\lambda_\infty(\pi_1\times\pi_2)$ where $q(\pi_1\times\pi_2)$ is the conductor appearing in the functional equation for $L(s,\pi_1\times\pi_2)$ and

\begin{equation*}
\lambda_\infty(\pi_1\times\pi_2)=\prod_{v=\infty}\prod_{i=1}^{n_1}\prod_{j=1}^{n_2}(1+|\mu_{\pi_1\times\pi_2}(v,i,j)|).
\end{equation*}

The definitions of $C(\pi)$ and $C(\pi_1\times\pi_2)$ were first made by Iwaniec and Sarnak in [7].  By the work of Bushnell and Henniart [1], $q(\pi_1\times\pi_2)\leq q(\pi_1)^{n_2}q(\pi_2)^{n_1}$.  It can also be shown that $\lambda_\infty(\pi_1\times\pi_2)\ll_{n_1,n_2} \lambda_\infty(\pi_1)^{n_2}\lambda_\infty(\pi_2)^{n_1}$.  Thus 

\begin{equation}\label{conductor}
C(\pi_1\times\pi_2)\ll_{n_1,n_2} C(\pi)^{n_2}C(\pi_2)^{n_1}.
\end{equation}

\bigskip

\qquad {\it Definition.}  Let $f(s,\pi)$ be a Dirichlet series associated with $\pi$ which converges absolutely to the right of 1.  We say that $f(s,\pi)$ {\it satisfies the convexity bound at $s=1$} if $f(s,\pi)=O_\epsilon (C(\pi)^\epsilon)$ for every $\epsilon>0$ and all ${\rm Re}(s)>1$.  When $f(s,\pi)$ has a functional equation and nice analytic properties which allow for an interpolation to the left of 1, we drop the reference to any particular point and say simply that $f(s,\pi)$ satisfies the {\it standard convexity bound}.

\bigskip

\qquad Our goal is to show that $L(s,\pi,|{\rm max}|^2)$ satisfies the convexity bound at $s=1$.  Uniform estimates in the conductor for a Dirichlet series are generally derived from a functional equation in which the conductor appears.  Unfortunately $L(s,\pi,|{\rm max}|^2)$ satisfies no such functional equation.  Proposition \ref{harken back} will allow us to obtain uniform estimates for $L(s,\pi,|{\rm max}|^2)$ from those for $L(s,\pi\times\widetilde\pi)$ and $L(s,\land^j\pi\times\widetilde{\land^j\pi})$ where $2\leq j\leq\lfloor n/2\rfloor$.

\bp\label{catherine} Let $\pi$ be a cuspidal representation of ${\rm GL}_n(\mathbb{A})$.  If $\land^j\pi$ is strongly automorphic isobaric for $2\leq j\leq\lfloor n/2\rfloor$, then $L(s,\pi,|{\rm max}|^2)=O(C(\pi)^A)$ on ${\rm Re}(s)>1$ for some $A>0$.\ep

\pf For $j\in\{ 2,\ldots,\lfloor n/2\rfloor\}$, let $\land^j\pi=\sigma_{j,1}\boxplus\cdots\boxplus\sigma_{j,\ell_j}$ be the decomposition of $\land^j\pi$ into an isobaric sum of cusp forms $\sigma_{j,i}$ on ${\rm GL}_{n_{j,i}}$.  Then 

\begin{equation}\label{product}
L(s,\land^j\pi\times\widetilde{\land^j\pi})=\prod_{1\leq i_1,i_2\leq\ell_j}L(s,\sigma_{j,i_1}\times\widetilde{\sigma_{j,i_2}}).
\end{equation}

The convergence of $L(s,\sigma_{j,i_1}\times\widetilde{\sigma_{j,i_2}})$ to the right of 1 along with its functional equation [8] imply, through the Phragmen-Lindelof convexity principle, that

\begin{equation*}
L(s,\sigma_{j,i_1}\times\widetilde{\sigma_{j,i_2}})=O(C(\sigma_{j,i_1}\times\widetilde{\sigma_{j,i_2}})^{B_{j,i_1,i_2}})
\end{equation*}

on ${\rm Re}(s)>1$ for some $B_{j,i_1,i_2}>0$.  By (\ref{conductor}) and (\ref{product}) we therefore have 

\begin{equation}\label{quick}
L(s,\land^j\pi\times\widetilde{\land^j\pi})=\sum_{\textfrak{n}}\lambda(\textfrak{n},\land^j\pi\times\widetilde{\land^j\pi})=O(C(\pi)^{B_j})
\end{equation}

on ${\rm Re}(s)>1$ where 

\begin{equation*}
B_j=\sum_{1\leq i_1,i_2\leq \ell_j}(n_{j,i_1}+n_{j,i_2})B_{j,i_1,i_2}.  
\end{equation*}

Similarly let $B_1>0$ be such that $L(s,\pi\times\tilde\pi)=O(C(\pi)^{B_1})$ on ${\rm Re}(s)>1$.  

\qquad Denote by $\Pi^{(j)}$ the strong exterior $j$ power lift of $\pi$.  Denote by $S_j$ the set of finite primes outside of which $\Pi^{(j)}_\textfrak{p}\simeq\land^j\pi_\textfrak{p}$.  Let $T$ be the union of all the $S_j$ and the set of primes at which $\pi$ is ramified.  As before, denote by $\textfrak{t}$ the product of all primes in $T$.  We note that ${\rm N}\textfrak{t}\ll C(\pi)$.

\qquad To bound $L_T(s,\pi,|{\rm max}|^2)$ (defined in (\ref{unramified L})) polynomially in $C(\pi)$ to the right of 1 we first note that for square-free ideals $\textfrak{n}$, $|\lambda(\textfrak{n},\pi,\land^j)|\leq 1+|\lambda(\textfrak{n},\pi,\land^j)|^2=1+\lambda(\textfrak{n},\land^j\pi\times\widetilde{\land^j\pi})$.  Secondly, Rudnick and Sarnak [17, Appendix] have shown that the coefficients of $L(s,\Pi\times\widetilde\Pi)$, where $\Pi$ is any isobaric form on ${\rm GL}_n$, are non-negative.  We may therefore remove the restriction of being square-free and relatively prime to $\textfrak{t}$.  Thus

\begin{equation*}
\sum_{(\textfrak{n},\textfrak{t})=1}^{\qquad\flat}|\lambda (\textfrak{n},\pi,\land^j)|{\rm N}\textfrak{n}^{-\sigma}\leq \sum_\textfrak{n}(1+\lambda (\textfrak{n},\pi,\land^j\times\widetilde{\land^j\pi})){\rm N}\textfrak{n}^{-\sigma}.
\end{equation*}

An appeal to Proposition \ref{harken back} and (\ref{quick}) gives $L_T(s,\pi,|{\rm max}|^2)=(C(\pi)^B)$ where $B=\sum_j B_j$ and $j$ runs through $1,\ldots ,\lfloor n/2\rfloor$. 

\qquad We now treat $L^T(s,\pi,|{\rm max}|^2)$ (defined in (\ref{ramified L})).  Let $\delta=\delta(n)=(n^2+1)^{-1}$.  Using the local bounds (\ref{LRS bounds}) we have

\begin{equation*}
L^\textfrak{p}(1,\pi,|{\rm max}|^2):=\sum_{r\geq 0}\underset{i}{\rm max}|\alpha_\pi (\textfrak{p},i)|^{2r}{\rm N}\textfrak{p}^{-r}\leq \sum_{r\geq 0}{\rm N}\textfrak{p}^{-2r\delta}=1+c{\rm N}\textfrak{p}^{-2\delta},
\end{equation*}

for some constant $c>0$ depending only on $\delta$.  Thus $L^\textfrak{p}(1,\pi,|{\rm max}|^2)\leq 1+{\rm N}\textfrak{p}^{-\delta}$ for primes $\textfrak{p}$ such that ${\rm N}\textfrak{p}\geq c^{\delta^{-1}}$.  We have

\begin{equation*}
L^T(1,\pi,|{\rm max}|^2)\underset{\delta}{\ll}\underset{\underset{{\rm N}\textfrak{p}\geq c^{\delta^{-1}}}{\textfrak{p}|\textfrak{t}}}{\prod}(1+{\rm N}\textfrak{p}^{-\delta})\leq \sum_{{\rm N}\textfrak{n}\leq {\rm N}\textfrak{t}}{\rm N}\textfrak{n}^{-\delta}\ll {\rm N}\textfrak{t}^{1-\delta}\underset{\delta}{\ll} C(\pi)^{1-\delta}.
\end{equation*}

With $A=B+1-\delta$, we have proved the proposition.\qed

\bigskip

\qquad Once we have polynomial control on $L(s,\pi,|{\rm max}|^2)$ to the right of 1 we can use a bootstrapping technique of Iwaniec [6] to whittle down the exponent to be as small as we like.  We now see the fruit of not having applied this technique straightaway to $L(s,\pi\times\pi)$, say, as in Molteni [16]: the complete multiplicativity of the coefficients $\lambda (\textfrak{n},\pi,|{\rm max}|^2)$ allows us to do without any further restriction on the size of the local roots.  That is, no improvement on the Luo-Rudnick-Sarnak bounds (\ref{LRS bounds}), already used in the proofs of Proposition \ref{harken back} and Proposition \ref{catherine}, will be necessary.  In fact, the following proposition could be applied to any Dirichlet series with completely multiplicative non-negative coefficients with polynomial control in the conductor to the right of 1.

\qquad We note that Propositions \ref{catherine} and \ref{barcelona} combine to give Theorem \ref{olivia}.

\bp\label{barcelona} Let $\pi$ be a cuspidal representation of ${\rm GL}_n(\mathbb{A})$.  Assume that the function $L(s,\pi,|{\rm max}|^2)$ converges on ${\rm Re}(s)>1$.  If there exists a constant $A>0$ such that $L(s,\pi,|{\rm max}|^2)=O(C(\pi)^A)$ on ${\rm Re}(s)> 1$ then $L(s,\pi,|{\rm max}|^2)$ satisfies the convexity bound at $s=1$.\ep

\pf For convenience, put $\lambda (\textfrak{n})=\lambda(\textfrak{n},\pi,|{\rm max}|^2)$ and $C=C(\pi)$.  Set $S(X)=\sum_{{\rm N}\textfrak{n}\leq X}\lambda(\textfrak{n})$.  Then the polynomial control and the positivity of the coefficients imply that for every $\sigma>1$

\begin{equation}\label{poly}
S(X)\leq X^\sigma \sum_{{\rm N}\textfrak{n}\leq X}\lambda(\textfrak{n}){\rm N}\textfrak{n}^{-\sigma}\leq X^\sigma \sum_\textfrak{n}\lambda(\textfrak{n}){\rm N}\textfrak{n}^{-\sigma}\ll C^AX^\sigma.
\end{equation}

\qquad By the complete multiplicativity of the $\lambda (\textfrak{n})$ we have

\begin{equation*}
S(X)^2 =\sum_{{\rm N}\textfrak{m},{\rm N}\textfrak{n}\leq X}\lambda (\textfrak{m})\lambda (\textfrak{n})=\sum_{{\rm N}\textfrak{m},{\rm N}\textfrak{n}\leq X}\lambda (\textfrak{mn})=\sum_{{\rm N}\textfrak{r}\leq X^2}\lambda(\textfrak{r})\tau(\textfrak{r}),
\end{equation*}

where $\tau (\textfrak{r})$ is the number of divisors of $\textfrak{r}$.  Applying the bound $\tau (\textfrak{r})\ll_\epsilon ({\rm N}\textfrak{r})^\epsilon$ and (\ref{poly}) we get $S(X)^2 \ll_\epsilon C^AX^{2+\epsilon}$.  Upon taking the square root we have $S(X)\ll_\epsilon C^{A/2}X^{1+\epsilon}$.  Iterating this step $M$ times, we obtain $S(X)\ll_{\epsilon, M} C^{A/2^M}X^{1+\epsilon}$.  For any $\epsilon>0$, we may take $M>(\log A-\log\epsilon)/\log 2$ to obtain

\begin{equation}\label{epsilon}
S(X)\underset{A, \epsilon}{\ll} C^\epsilon X^{1+\epsilon}.
\end{equation}

\qquad Let ${\rm N}\textfrak{n}\sim M$ denote the diadic interval $M\leq {\rm N}\textfrak{n}< 2M$.  Using (\ref{epsilon}) along with the positivity of the coefficients we conclude that, for any $\epsilon>0$ and $\sigma\geq1+2\epsilon$, $L(\sigma,\pi,|{\rm max}|^k)$ is

\begin{equation*}
\underset{k\geq 0}{\underset{M=2^k}{\sum}}\sum_{{\rm N}\textfrak{n}\sim M} \lambda(\textfrak{n}){\rm N}\textfrak{n}^{-\sigma}\leq\underset{k\geq 0}{\underset{M=2^k}{\sum}}M^{-\sigma}S(2M)\underset{A, \epsilon}{\ll} C^\epsilon\sum_{k\geq 0}2^{-k\epsilon}\underset{A, \epsilon}{\ll} C^\epsilon.
\end{equation*}

This finishes the proof.\qed

\bigskip

\qquad Theorem \ref{olivia} has now been proven.  Note that even for large $n$ the hypothesis that all $\land^j\pi$, $2\leq j\leq\lfloor n/2\rfloor$, be automorphic is not necessarily stronger than the condition $\alpha_\pi (\textfrak{p},i)\ll {\rm N}\textfrak{p}^{1/4}$ on the local roots.  Corollary \ref{over with} states that only when these exterior power lifts are combined with the {\it symmetric square} lift do the hypothetical bounds (\ref{molt}) follow.  And yet the strength of the conclusion of Theorem \ref{olivia} is much stronger than the convexity bound for only $L(s,\pi\times\pi)$.

\subsection{\rm Proof of Corollary \ref{nice}}\label{cor-proof}

\qquad When Theorem \ref{olivia} is combined with the (strong) automorphy of $\land^2\pi$ for $\pi$ on ${\rm GL}_4$, a fact proved in [9, Theorem 5.3.1], we obtain the following corollary.

\bigskip

\qquad {\sc Corollary \ref{nice}.}  {\it Let $\pi_i$ be cuspidal representations of ${\rm GL}_{n_i}(\mathbb{A})$ where $n_i\leq 4$ for $i=1, 2$.  Then $L(s,\pi_1\times\pi_2)$, as well as $L(s,\pi_i,\land^2)$ and $L(s,\pi_i,{\rm sym}^2)$ for $i=1, 2$, satisfy the standard convexity bound.}

\bigskip

\pf For any integer $r\geq 0$ and prime ideal $\textfrak{p}$ we have

\begin{equation*}
\lambda(\textfrak{p}^r,\pi\times\tilde\pi)\leq \mathcal{N}(r,n)\  \underset{i}{\rm max}|\alpha_\pi (\textfrak{p},i)|^{2r},
\end{equation*}

where $\mathcal{N}(r,n)$ is the number of monomials in $n$ variables of degree $r$.  The same bound holds for $\lambda(\textfrak{p}^r,\pi,\land^2)$ and $\lambda(\textfrak{p}^r,\pi,{\rm sym}^2)$.  We can compute that $\mathcal{N}(r,n)\ll r^{A(n)}$ for some $A(n)>0$.  Theorem \ref{olivia} therefore gives the convexity bound for each of $L(s,\pi_i\times\widetilde\pi_i)$, $L(s,\pi_i,\land^2)$, and $L(s,\pi_i, {\rm sym}^2)$.  

\qquad We deduce the convexity bound for $L(s,\pi_1\times\pi_2)$ from that of $L(s,\pi_i\times\widetilde\pi_i)$ for $i=1,2$.  To do so, we avail ourselves of the notation and terminology of Macdonald's treatise [15].  For integers $n\geq 1$ and $r\geq 0$, let $\mathcal{P}_n(r)$ be the set of partitions of $r$ of length no greater than $n$.  Let $s_\lambda$ denote the Schur function associated to a partition $\lambda$.  Assume $n_1\geq n_2$.  If $n_2 < n_1$ then define $\alpha_{n_2}(\textfrak{p},i)=0$ for all $n_2<i\leq n_1$.  Put $n=n_1$ and set $\alpha_{\pi_i} (\textfrak{p})=(\alpha_{\pi_i} (\textfrak{p},1),\ldots ,\alpha_{\pi_i} (\textfrak{p},n))$.  For primes $\textfrak{p}$ unramified for both $\pi_1$ and $\pi_2$ and integers $r\geq 0$ the coefficients $\lambda (\textfrak{p}^r,\pi_1\times\pi_2)$ are defined by

\begin{equation*}
\prod_{1\leq i_1,i_2\leq n}(1-\alpha_{\pi_1}(\textfrak{p},i_1)\alpha_{\pi_2}(\textfrak{p},i_2){\rm N}\textfrak{p}^{-s})^{-1}=\sum_{r\geq 0}\lambda(\textfrak{p}^r,\pi_1\times\pi_2){\rm N}\textfrak{p}^{-rs}.
\end{equation*}

It is then a standard identity in the theory of symmetric functions that

\begin{equation*}
\lambda(\textfrak{p}^r,\pi_1\times\pi_2)=\sum_{\lambda\in\mathcal{P}_n(r)}s_\lambda(\alpha_{\pi_1}(\textfrak{p}))s_\lambda (\alpha_{\pi_2}(\textfrak{p})).
\end{equation*}

Applying the Cauchy-Schwartz inequality to this, we obtain $|\lambda (\textfrak{p}^r,\pi_1\times\pi_2)|\leq \lambda(\textfrak{p}^r,\pi_1\times\widetilde\pi_1)^{1/2}\lambda (\textfrak{p}^r,\pi_2\times\widetilde\pi_2)^{1/2}$.  A similiar inequality can be proven for ramified primes.  When extended to all $\textfrak{n}$ this gives

\begin{equation*}
\left|\sum_\textfrak{n}\frac{\lambda (\textfrak{n},\pi_1\times\pi_2)}{\textfrak{n}^s}\right|\leq \sum_{\textfrak{n}}\frac{\lambda(\textfrak{n},\pi_1\times\widetilde\pi_1)^{1/2}}{\textfrak{n}^{\sigma/2}}\frac{\lambda(\textfrak{n},\pi_2\times\widetilde\pi_2)^{1/2}}{\textfrak{n}^{\sigma/2}}\leq L(s,\pi_1\times\widetilde\pi_1)L(s,\pi_2\times\widetilde\pi_2),
\end{equation*}

the last inequality again by Cauchy-Schwartz.  The corollary immediately follows.\qed

\begin{center}
{\sc References}
\end{center}

[1] C.J. Bushnell and G. Henniart, {\it An upper bound on conductors for pairs,} J. Number Theory 65 (1997), no. 2, 183-196.

\medskip

[2] W. Duke and E. Kowalski, {\it A problem of Linnik for elliptic curves and mean-value estimates for automorphic representations,}  With an appendix by Dinakar Ramakrishnan.  Invent. Math. {\bf 139} (2000), no. 1, 1--39.

\medskip

[3] S. Gelbart and H. Jacquet, {\it A relation between automorphic representations of GL(2) and GL(3),} Ann. Sci. Ecole Norm. Sup. (4) 11 (1978), no. 4, 471-542.

\medskip

[4] M. Harris and R. Taylor, {\it On the geometry and cohomology of some simple Shimura varieties,} Annals of Mathematics Studies, 151, Princeton University Press, 2001.

\medskip

[5] G. Henniart, {\it Une preuve simple des conjectures de Langlands pour GL(n) sur un corps $p$-adique,} Inv. Math. 139 (2000), 439-455.

\medskip

[6]  H. Iwaniec, {\it Small Eigenvalues of Laplacian for $\Gamma_0 (N)$,} Acta Arith. {\bf 56} (1990), no.1, 65-82.

\medskip

[7] H. Iwaniec and P. Sarnak, {\it Perspectives on the analytic theory of $L$-functions,} GAFA 2000 (Tel Aviv, 1999). Geom. Funct. Anal.  {\bf 2000},  Special Volume, Part II, 705--741.

\medskip

[8] H. Jacquet, I. Piatetskii-Shapiro, and J. Shalika, {\it Rankin-Selberg convolutions,}  Amer. J. Math. {\bf 105} (1983), no. 2, 367--464.

\medskip

[9]  H. Kim, {\it Functoriality for the exterior square of ${\rm GL}_4$ and the symmetric fourth of ${\rm GL}_2$,}  J. Amer. Math. Soc. {\bf 16} (2003), no.1, 139-183.

\medskip

[10] E. Kowalski and P. Michel, {\it Zeros of families of automorphic $L$-functions close to 1}, Pac. Journ. of Math., Vol. 207, no. 2, 2002, 411-431.

\medskip

[11] H. Kim and F. Shahidi, {\it Cuspidality of symmetric powers with applications,} Duke Math. J., 112 (1) 2002, 177-197.

\medskip

[12] R.P. Langlands, {\it On the classification of irreducible representations of real algebraic groups,} in Representation Theory and Harmonic Analysis on Semisimple Lie groups (P.J. Sally, Jr. and D.A. Vogan, ed.), Mathematical Surveys and Monographs, vol. 31, AMS, 1989, pp. 101-170.

\medskip

[13] W. Luo, {\it Values of symmetric square $L$-functions at 1}, J. Reine Angew. Math., {\bf 506} (1999), 215-235.

\medskip

[14] W. Luo, Z. Rudnick and P. Sarnak, {\it On the generalized Ramanujan conjecture for ${\rm GL}(n)$,} in Automorphic forms, automorphic representations, and arithmetic, Proc. Sympos. Pure Math., vol. 66, Part 2, Amer. Math. Soc., Providence, RI, 1999, pp. 301-310.

\medskip

[15]  Macdonald, I. G. Symmetric functions and Hall polynomials. Second edition. With contributions by A. Zelevinsky. Oxford Mathematical Monographs. Oxford Science Publications. {\it The Clarendon Press, Oxford University Press, New York,} 1995. x+475 pp.

\medskip

[16] G. Molteni, {\it Upper and lower bounds at $s=1$ for certain Dirichlet series with Euler product,} Duke Math J., Vol. 111, No. 1 (2000), 133-158.

\medskip

[17] Z. Rudnick and P. Sarnak, {\it Zeros of principal $L$-functions and random matrix theory.} A celebration of John F. Nash, Jr.,  Duke Math. J. {\bf 81} (1996), no.2, 269-322.

\bigskip

\qquad {\sc I3M, UMR CNRS 5149, Universit\'e Montpellier II CC 051, 34095 Montpellier Cedex 05, France}

\qquad {\it E-mail address:} {\tt brumley@math.univ-montp2.fr}

\end{document}